\documentclass{amsart}

\usepackage{amssymb}

\newtheorem{theorem}{Theorem}[section]
\newtheorem{Prop}[theorem]{Proposition}
\newtheorem{Thm}[theorem]{Theorem}
\newtheorem{Lem}[theorem]{Lemma}

\theoremstyle{definition}

\newtheorem{Def}[theorem]{Definition}

\theoremstyle{remark}

\numberwithin{equation}{section}

\newcommand{\R}{{\mathbb R}}
\newcommand{\C}{{\mathbb C}}
\newcommand{\Z}{{\mathbb Z}}

\newcommand{\LG}{{\mathfrak g}}

\newcommand{\LN}{{\mathfrak n}}

\newcommand{\LV}{{\mathfrak v}}
\newcommand{\LZ}{{\mathfrak z}}

\newcommand{\SO}{\mathrm{SO}} 
\newcommand{\OO}{\mathrm{O}} 
\newcommand{\GL}{\mathrm{GL}} 
\newcommand{\CO}{\mathrm{CO}} 
 
\newcommand{\SL}{\mathrm{SL}}

\newcommand{\im}{\mathrm{Im} \ }
\newcommand{\re}{\mathrm{Re} \ }
\newcommand{\Aut}{\mbox{\rm Aut}}
\newcommand{\Cliff}{\mathrm{Cliff}}
\newcommand{\rank}{\mathrm{rank} \ }
\newcommand{\End}{\mbox{\rm End}}
\newcommand{\Hom}{\mathrm{Hom}}

\newcommand{\Pin}{\mathrm{Pin}}

\newcommand{\id}{\mbox{\rm id}}

\newcommand{\htype}{{\it H}-type } 
\newcommand{\inner}[2]{\langle #1 , #2 \rangle} 

\begin{document}

\title{Lie groups locally isomorphic to generalized Heisenberg groups}

\author{Hiroshi Tamaru}
\address{Department of Mathematics, Hiroshima University, Higashi-Hiroshima 739-8526, Japan}
\email{tamaru@math.sci.hiroshima-u.ac.jp}
\thanks{The first author was supported in part by 
Grant-in-Aid for Young Scientists (B) 17740039, 
The Ministry of Education, Culture, Sports, Science and Technology, Japan.}

\author{Hisashi Yoshida}
\address{
Graduate School of Mathematics, 
Nagoya University, 
Chikusa-ku, Nagoya 464-8602, 
Japan
}
\email{r060141m@mbox.nagoya-u.ac.jp}


\begin{abstract}
We classify connected Lie groups which are locally isomorphic to 
generalized Heisenberg groups. 
For a given generalized Heisenberg group $N$, 
there is a one-to-one correspondence between 
the set of isomorphism classes of connected Lie groups 
which are locally isomorphic to $N$ 
and a union of certain quotients of noncompact Riemannian symmetric spaces. 
\end{abstract}

\maketitle


\section{Introduction}

Generalized Heisenberg groups were introduced by Kaplan (\cite{Kaplan}), 
and have been studied in many fields in mathematics. 
Especially, generalized Heisenberg groups, endowed with the left-invariant metric, 
provide beautiful examples in geometry. 
They provide many examples of symmetric-like Riemannian manifolds, 
that is, Riemannian manifolds which share some properties with 
Riemannian symmetric spaces 
(see \cite{BTV} and the references). 
Compact quotients of some generalized Heisenberg groups 
provide isospectral, but nonisometric manifolds 
(see Gordon \cite{Gordon} and the references). 
Generalized Heisenberg groups also provide remarkable examples of solvmanifolds, 
by taking certain solvable extension. 
The constructed solvmanifolds are now called Damek-Ricci spaces, 
after Damek-Ricci (\cite{DR}) showed that they are harmonic spaces 
(see also \cite{BTV}). 

The aim of this paper is to classify connected Lie groups 
which are locally isomorphic to generalized Heisenberg groups. 
They provide interesting examples of nilmanifolds, 
that is, Riemannian manifolds on which nilpotent Lie groups act 
isometrically and transitively. 
They might be a good prototype of the study of 
non simply-connected and noncompact nilmanifolds. 
This is one of our motivation. 

Other motivation comes from the interesting structure of the moduli space, 
the set of isomorphism classes of connected Lie groups 
locally isomorphic to a generalized Heisenberg group. 
If one thought about abelian Lie groups or classical Lie groups, 
the moduli space seems to be a finite set. 
But, in the case of a generalized Heisenberg group $N$, 
there is a continuous family of connected Lie groups 
which are locally isomorphic to given $N$ 
(except for the classical Heisenberg groups). 
Our main theorem states that the moduli space is bijective to 
\[ 
\coprod_{r=0}^{m} \ ( \SL_r(\Z) \backslash \SL_r(\R) / \SO_r(\R) ) , 
\] 
where $m$ is the dimension of the center $Z(N)$. 
It seems to be surprising that the moduli spaces depend only on $\dim Z(N)$. 

This paper is organized as follows. 
In Section \ref{section_auto}, 
we study the structure of the automorphism group of a metric Lie algebra 
and of a two-step nilpotent metric Lie algebra. 
We recall the definition of a generalized Heisenberg group, 
and describe the action of the automorphism group on the center, 
in Section \ref{section_htype}. 
In particular, this action is equivalent to the action of 
the conformal group $\CO_m(\R)$. 
Our main theorem is proved in Section \ref{section_main}. 
For the proof, 
we have to classify lattices in $\R^m$ up to $\CO_m(\R)$-conjugation. 
In case of $\dim Z(N) =2$, 
the component $\SL_2(\Z) \backslash \SL_2(\R) / \SO_2(\R)$ 
is the quotient of the real hyperbolic plane $\R H^2$ by the 
modular group $\SL_2(\Z)$, which have been well studied. 
Using the structure theorem of $\SL_2(\Z) \backslash \R H^2$, 
we give explicit descriptions 
of connected Lie groups which are locally isomorphic to $N$ in this case, 
in Section \ref{section_two}. 

The authors would like to thank 
Professor Makoto Matsumoto, Professor Nobuo Tsuzuki, 
Dr.\ Akira Ishii, Dr.\ Nobuyoshi Takahashi, Dr.\ Masao Tsuzuki 
and Dr.\ Takuya Yamauchi 
for their useful advices and kind encouragements. 
The authors are also grateful to 
Mr.\ Tadashi Kashiwa and Mr.\ Hironao Kato for useful discussions. 
%
%

\section{Preliminaries on automorphism groups}
\label{section_auto}

For the classification of connected Lie groups locally isomorphic to a given 
simply-connected Lie group $G$, 
we need to know the action of the automorphism group on the center of $G$. 
In this section we study the structure of the automorphism group 
of a two-step nilpotent metric Lie algebra. 

We start with an arbitrary metric Lie algebra $(\LG, \inner{}{})$, 
a Lie algebra with an inner product. 
Let $\LZ$ be the center of $\LG$ and 
$\LV$ the orthogonal complement of $\LZ$ in $\LG$. 
Since the automorphism group $\Aut(\LG)$ preserves $\LZ$, 
we first study the structure of 
\begin{equation*}
N_{\GL(\LG)}(\LZ) := \{ \sigma \in \GL(\LG); \ \sigma(\LZ) \subset \LZ \} . 
\end{equation*}
Let us define an additive group 
\begin{equation*}
\Hom(\LV,\LZ) := 
\{ \beta : \LV \rightarrow \LZ : \mbox{linear} \} . 
\end{equation*}
Each element $\beta \in \Hom(\LV,\LZ)$ acts on $\LG$ by $T_{\beta}$ 
in the natural way, that is, 
\begin{equation*}
T_{\beta}(z) :=  z, \ T_{\beta}(v) := \beta(v) + v, \quad 
\mbox{for $z \in \LZ$ and $v \in \LV$}. 
\end{equation*} 

\begin{Lem}
\label{lemma1}
$N_{\GL(\LG)}(\LZ) = \Hom(\LV,\LZ) \rtimes (\GL(\LZ) \times \GL(\LV))$. 
\end{Lem}

\begin{proof}
Let $(h,g) \in \GL(\LZ) \times \GL(\LV)$ and $\beta \in \Hom(\LV,\LZ)$. 
It is easy to see that 
$(h,g)^{-1} \circ T_{\beta} \circ (h,g) = T_{h^{-1} \circ \beta \circ g}$. 
Thus $\GL(\LZ) \times \GL(\LV)$ normalizes $\Hom(\LV,\LZ)$ 
and the inclusion "$\supset$" is clear. 
To show the converse, let $T \in N_{\GL(\LG)}(\LZ)$. 
Since $T$ preserves $\LZ$, 
there exist $h \in \GL(\LZ)$, $g \in \GL(\LV)$ and $\beta \in \Hom(\LV,\LZ)$ 
such that 
\begin{equation*}
T(z) = h(z), \quad T(v) = \beta(v) + g(v). 
\end{equation*}
It concludes that 
$T = T_{\beta \circ g^{-1}} \circ (h,g) \in 
\Hom(\LV,\LZ) \rtimes (\GL(\LZ) \times \GL(\LV))$. 
\end{proof}

Let us now consider two-step nilpotent metric Lie algebras $(\LN, \inner{}{})$. 
In this case $\Aut(\LN)$ can be decomposed into two groups, 
according to Lemma \ref{lemma1}. 
We need 
\begin{equation*}
\Aut_{\LV}(\LN) := \Aut(\LN) \cap (\GL(\LZ) \times \GL(\LV)) . 
\end{equation*}

\begin{Lem}
\label{dec}
For a two-step nilpotent metric Lie algebra $(\LN, \langle, \rangle)$, 
we have 
\begin{enumerate}
\item 
$\Hom(\LV,\LZ) \subset \Aut(\LN)$, and therefore 
$\Aut(\LN) = \Hom(\LV,\LZ) \rtimes \Aut_{\LV}(\LN)$, and 
\item
$\Aut_{\LV}(\LN) = 
\{(h,g) \in \GL(\LZ) \times \GL(\LV); \ 
[gu,gv] = h [u,v] \  (\forall u,v \in \LV) \}$. 
\end{enumerate}
\end{Lem}

\begin{proof}
Since $\LN$ is of two-step nilpotent, 
one can immediately see that $[\LV, \LV] \subset \LZ$. 
This leads easily that $\Hom(\LV,\LZ) \subset \Aut(\LN)$. 
Lemma \ref{lemma1} states that 
$\Aut(\LN) \subset \Hom(\LV,\LZ) \rtimes (\GL(\LZ) \times \GL(\LV))$. 
We can conclude (1) by the following group-theoretic property: \ 
Let $K \subset H_1 \rtimes H_2$ and $H_1 \subset K$, 
then $K = H_1 \rtimes (H_2 \cap K)$ holds. 
The claim (2) is a easy consequence of 
the assumption that $\LN$ is of two-step nilpotent. 
\end{proof}

We note that Lemma \ref{dec} is obtained by Saal (\cite[Proposition 2.3]{Saal}) 
for generalized Heisenberg algebras. 

\section{Automorphism groups of \htype groups}
\label{section_htype}

In this section, 
we recall the definition of a generalized Heisenberg algebra 
$(\LN, \langle, \rangle)$ 
and describe the action of the automorphism group $\Aut(\LN)$ 
on the center of $\LN$. 

\begin{Def} 
A two-step nilpotent metric Lie algebra $(\LN, \langle, \rangle)$ is called 
a {\it generalized Heisenberg algebra} or an {\it \htype algebra} if 
\[ 
J_z^2 = -\langle z,z \rangle \id_{\LV} \quad (\forall z \in \LZ), 
\] 
where the operator $J : \LZ \rightarrow \End(\LV)$ is defined by 
\[ 
\langle J_zv,u \rangle = \langle z,[v,u] \rangle \quad 
\mbox{for $z \in \LZ$, $u,v \in \LV$}. 
\] 
The connected and simply-connected Lie group $N$ with Lie algebra $\LN$, 
endowed with the induced left-invariant metric, is called a 
{\it generalized Heisenberg group} or an {\it \htype group}. 
\end{Def}

We refer \cite{Kaplan} and \cite{BTV} for \htype groups and algebras. 

We would like to know the action of the automorphism group $\Aut(N)$ 
on the center $Z(N)$ of an \htype group $N$. 
Since $N$ is simply-connected and the exponential map 
$\exp : \LN \rightarrow N$ is a diffeomorphism, 
it is sufficient to investigate the action of the automorphism group $\Aut(\LN)$ 
on the center $\LZ$ of $\LN$. 
Note that the following proposition can be obtained as a corollary of 
the description of the automorphism group $\Aut(\LN)$ by Saal (\cite{Saal}). 
But, we will give a proof here, 
since we determine the action of $\Aut(\LN)$ on $\LZ$ directly 
and hence the arguments become slightly simpler. 
Denote the conformal group by 
\[ 
\CO(\LZ) := \{ dg \in \GL(\LZ) \mid d \in \R^{\times}, \  g \in \OO(\LZ) \} . 
\] 

\begin{Prop}
\label{D-RO}
For an \htype algebra $(\LN, \langle, \rangle)$, 
the action of $\Aut(\LN)$ on the center $\LZ$ 
is equivalent to the action of the conformal group $\CO(\LZ)$ on $\LZ$. 
\end{Prop}

\begin{proof}
First of all, Lemma \ref{dec} states that 
$\Aut(\LN) = \Hom(\LV,\LZ) \rtimes \Aut_{\LV}(\LN)$. 
By definition, $\Hom(\LV,\LZ)$ acts trivially on $\LZ$, 
therefore the actions of $\Aut(\LN)$ and $\Aut_{\LV}(\LN)$ coincide on $\LZ$. 
Hence, the group we would like to know is nothing but 
\[ 
D:= \{h \in \GL(\LZ); \ \exists g \in \GL(\LV) : 
(h,g) \in \mathrm{Aut}_\LV(\LN) \} . 
\] 
For the proof of $D = \CO(\LZ)$, we need the following groups: 
\begin{itemize}
\item 
Let $\Cliff(\LN)$ denote the subgroup of $\Aut(\LN)$ generated by 
$\{ (-|z|^2 \rho_z, J_z); \ z \in \LZ -\{0\} \}$, 
where $\rho_z$ is the reflection in $\LZ$ with respect to $(\R z)^{\perp}$. 
\item
Let $\Pin(\LN)$ be the subgroup of $\Cliff(\LN)$ generated by 
$\{(-\rho_z,J_z); \ |z| = 1\}$. 
\end{itemize}
Note that $(-|z|^2 \rho_z, J_z) \in \GL(\LZ) \times \GL(\LV)$ 
defines an automorphism of $\LN$. 
Let $l := \dim \LV$, $m := \dim \LZ$, and 
$\{ z_1, \ldots, z_m \}$ be an orthonormal basis of $\LZ$. 

Claim 1: 
the action of $\Cliff(\LN)$ on $\LZ$ coincides with that of $\CO(\LZ)$. 
In case of $m$ is even, 
the action of $\Pin(\LN)$ on $\LZ$ coincides with that of $\OO(\LZ)$, 
since $\det(- \rho_z) = -1$. 
Therefore $\Cliff(\LN)$ acts on $\LZ$ as $\CO(\LZ)$. 
In case of $m$ is odd, 
the action of $\Pin(\LN)$ on $\LZ$ coincides with that of $\SO(\LZ)$, 
since $\det(- \rho_z) = 1$. 
Thus $\Cliff(\LN)$ acts on $\LZ$ as $\R^{\times} \cdot \SO(\LZ)$, 
but $m$ is odd, which implies 
$\R^{\times} \cdot \SO(\LZ) = \R^{\times} \cdot \OO(\LZ) = \CO(\LZ)$. 

Claim 2: 
$D \supset \CO(\LZ)$. 
This is a direct consequence of Claim 1. 

Claim 3: 
$D \subset \CO(\LZ)$. 
Let us take $h \in D$. 
By definition, 
there exists $g \in \GL(\LV)$ such that $(h,g) \in \Aut_{\LV}(\LN)$. 
Elementary linear algebra leads that 
\[ 
\exists k_1, k_2 \in O(\LZ) \ \mbox{such that} \ 
k_1 h k_2 = \mathrm{diag}(d_1,\ldots,d_m) \ \mbox{and} \ d_i > 0. 
\] 
Put $h_0 := \mathrm{diag}(d_1,\ldots,d_m)$, which is the diagonal matrices 
(note that we consider matrices representation with respect to 
the basis $\{ z_1, \ldots, z_m\}$). 
Since Claim 2 implies that 
$\Cliff(\LN)$ has a subgroup which acts on $\LZ$ as $O(\LZ)$, 
there exists $g_1, g_2 \in \GL(\LV)$ such that 
$(k_1,g_1)$, $(k_2,g_2) \in \Cliff(\LN) \subset \Aut_{\LV}(\LN)$. 
Let $g_0 := g_1gg_2$. One has 
\[ 
(h_0,g_0) = (k_1,g_1)(h,g)(k_2,g_2) \in \Aut_{\LV}(\LN) . 
\] 
The condition for $(h_0,g_0)$ to be an automorphism is 
\[ 
[g_0u,g_0v] = h_0 [u,v] \quad (\forall u,v \in \LV) . 
\] 
By taking the inner product with $z_i$ ($i = 1, \ldots, m$), we have 
\[ 
\langle z_i, [g_0u,g_0v] \rangle = \langle z_i, h_0 [u,v] \rangle \quad 
(\forall u,v \in \LV). 
\] 
This leads that 
${}^tg_0J_{z_i}g_0 = d_iJ_{z_i}$. 
By taking the determinant of each sides, one has $(\det g_0)^2 = d_i^l$. 
This means that $d_i$ is independent of $i$ (recall that $d_i > 0$). 
Therefore $h_0 = d \cdot \id_{\LZ}$, 
where $d := d_1 = d_2 = \cdots = d_m$. 
We conclude that $h = dk_1^{-1}k_2^{-1} \in \CO(\LZ)$, 
which completes the proof. 
\end{proof}

\section{Main Theorem} 
\label{section_main}

In this section, we classify connected Lie groups 
which are locally isomorphic to \htype groups $N$. 
We refer \cite{Hel}, \cite{Sugiura} for the theory of Lie groups and 
homogeneous spaces. 

Let $G$ be a connected and simply-connected Lie group, and denote by 
\[ 
LI(G) := \{ G'; \ \mbox{$G'$ is a connected Lie group and locally isomorphic to $G$} \} 
/ \mbox{isom.}
\]
An arbitrary connected Lie group can be represented by $G / \Gamma$, 
where $G$ is a connected and simply-connected Lie group 
and $\Gamma$ is a discreate subgroup contained in the center $Z(G)$ of $G$. 
Furthermore, $G / \Gamma_1$ and $G / \Gamma_2$ are isomorphic to each other 
if and only if $\Gamma_1$ and $\Gamma_2$ are conjugate by $\Aut(G)$. 
Therefore, one has 
\[ 
LI(G) \cong \{ \Gamma \subset Z(G); \ 
\mbox{$\Gamma$ is a discreate subgroup} \} / \Aut(G). 
\] 
We now study \htype groups $N$. 
Denote the Lie algebra of $N$ by $\LN = \LZ \oplus \LV$, 
where $\LZ$ is the Lie algebra of $Z(N)$. 
Via a Lie group isomorphism $\exp : \LZ \to Z(N)$, 
discreate subgroups in $Z(N)$ corresponds to 
lattices in $\LZ \cong \R^m$, where $m = \dim \LZ$. 
Note that $\Aut(N) \cong \Aut(\LN)$, since $N$ is simply-connected. 
Therefore we have that 
\[ 
LI(N) \cong \{ \Lambda \subset \LZ; \ \mbox{$\Lambda$ is a lattice} \} / \Aut(\LN). 
\] 
Hence, Proposition \ref{D-RO} leads that 
\[ 
LI(N) \cong \{ \Lambda \subset \LZ; \ \mbox{$\Lambda$ is a lattice} \} / \CO(\LZ) . 
\] 

\begin{Thm}
\label{maintheorem}
Let $N$ be an \htype group with $m$-dimensional center. 
The moduli space $LI(N)$ of isomorphism classes of connected Lie groups 
which are locally isomorphic to $N$ satisfies 
\[ 
LI(N) \cong \coprod_{r=0}^{m} \ ( \SL_r(\Z) \backslash \SL_r(\R) / \SO_r(\R) ) , 
\] 
where we define $\SL_0(\R) \backslash \SL_0(\Z) / \SO_0(\R) = \{ 0 \}$. 
\end{Thm}

\begin{proof}
We fix an orthonormal basis $\{ z_1, \ldots, z_m \}$ of $\LZ$. 
Then one can identify $\LZ = \R^m$, and thus 
\[ 
LI(N) = \{ \Lambda \subset \R^m; \ \mbox{$\Lambda$ is a lattice} \} / \CO_m(\R) . 
\] 
Here we define 
\[ 
LI_r(N) := \{ \Lambda \subset \R^m; \ \mbox{$\Lambda$ is a lattice of rank $r$} \} / \CO_m(\R) . 
\] 
On the other hand, 
the natural inclusion map $\SL_r(\R) \rightarrow \GL_r(\R)$ induces 
\[ 
\SL_r(\Z) \backslash \SL_r(\R) / \SO_r(\R) \cong 
\GL_r(\Z) \backslash \GL_r(\R) / \CO_r(\R) , 
\]
since $\CO_r(\R) \cap \SL_r(\R) = \SO_r(\R)$ 
and $\GL_r(\Z) \cap \SL_r(\R) = \SL_r(\Z)$. 
We will show the theorem by constructing a bijective map 
from $\GL_r(\Z) \backslash \GL_r(\R) / \CO_r(\R)$ to $LI_r(N)$. 
First of all, we define 
\[ 
\varphi : \GL_r(\R) \rightarrow 
\{ \Lambda \subset \R^m; \ \rank \Lambda = r \} : 
\left[ 
\begin{array}{c} 
u_1 \\ \vdots \\ u_r 
\end{array} 
\right] 
\mapsto \Z \tilde{u}_1 + \cdots + \Z \tilde{u}_r, 
\] 
where $u_i \in M_{1,r}(\R)$, 
$\tilde{u}_i := (u_i, 0, \ldots, 0) \in M_{1,m}(\R)$. 
Denote by $\tilde{\varphi}$ the induced map, 
\[ 
\tilde{\varphi} : 
\GL_r(\Z) \backslash \GL_r(\R) / \CO_r(\R) \rightarrow LI_r(N) : 
u \mapsto [\varphi(u)], 
\] 
where $[\varphi(u)]$ is the $\CO_m(\R)$-conjugate class of the lattice 
$\varphi(u)$. 

Claim 1: $\tilde \varphi$ is well-defined. 
Let $u \in \GL_r(\R)$, $C \in \GL_r(\Z)$ and $g \in \CO_r(\R)$. 
One has $[\varphi(Cug)] = [\varphi(Cu)]$, 
since $\CO_r(\R)$ is a subgroup of $\CO_m(\R)$ naturally. 
It can also be easily checked that $\varphi(u) = \varphi(v)$ 
(that is, they are the same lattice) if and only if 
$u = Cv$ for some $C \in \GL_r(\Z)$. 
Hence one has $\varphi(Cu) = \varphi(u)$. 
We conclude that $[\varphi(Cug)] = [\varphi(u)]$. 

Claim 2: $\tilde \varphi$ is surjective. 
Let $\Lambda = \Z u'_1 + \cdots + \Z u'_r$ be a lattice 
so that $[\Lambda] \in LI_r(N)$. 
Then there exists $g \in \CO_m(\R)$ such that
$u'_ig = (u_i, 0, \ldots, 0)$ (for every $i=1, \ldots, r$). 
Since $\rank \Lambda =r$, one has 
\[ 
u := 
\left[ 
\begin{array}{c} 
u_1 \\ \vdots \\ u_r 
\end{array} 
\right] 
\in \GL_r(\R)
\] 
and $\varphi(u) g^{-1} = \Lambda$. 
We conclude that $[\varphi(u)] = [\Lambda]$. 

Claim 3: $\tilde \varphi$ is injective. 
Let $u,v \in \GL_r(\R)$, and 
assume that 
$[\varphi(u)] = [\varphi(v)]$. 
Then there exists $g \in \CO_m(\R)$ such that 
$\varphi(u) = \varphi(v) g$. 
Since all the basis of these lattices are of the form $(u_i, 0, \ldots, 0)$, 
one can see that there exists $g' \in \CO_r(\R)$ such that 
$\varphi(u) = \varphi(v)g' = \varphi (vg')$. 
They are the same lattice, which implies that 
$u = Cvg'$ for some $C \in \GL_r(\Z)$. 
We conclude that $u$ and $v$ are in the same equivalent class. 
\end{proof}

\section{The case of two dimensional center} 
\label{section_two}

In this section we explicitly describe the classification of 
Lie groups which are locally isomorphic to \htype groups 
with two dimensional center. 

It is known that an \htype group $N$ 
with two dimensional center is the complex Heisenberg group, that is, 
\[ 
N = \left\{ \left[ 
\begin{array}{ccccc}
1 & \alpha_1 & \cdots & \alpha_k & \gamma \\ 
& 1 & & & \beta_k \\ 
& & \ddots & & \vdots \\ 
& & & 1 & \beta_1 \\ 
& & & & 1 
\end{array}
\right] \ ; \ \alpha_i, \beta_i, \gamma \in \C \right\} . 
\] 
We denote the element in the center $\LZ$ of the Lie algebra $\LN$, 
for $\gamma \in \C$, by 
\[
u(\gamma) := 
\left[ 
\begin{array}{cccc}
0 & \cdots & 0 & \gamma \\ 
& \ddots & & 0 \\ 
& & \ddots & \vdots \\ 
& & & 0 
\end{array}
\right] . 
\]
Obviously $\{ u(1), u(i) \}$ is an orthonormal basis of $\LZ$. 
Define 
\[ 
D_1 := \{ \rho \in \C; \ \im \rho >0, \ 
|\rho| \leq 1, \ |\rho -1| \geq 1, \ |\rho +1| \geq 1 \} . 
\]

\begin{Thm}
Let $N$ be a complex Heisenberg group. 
A connected (real) Lie group which is locally isomorphic to $N$ is 
isomorphic to either 
\begin{enumerate}
\item
$N$, 
\item
$N/ \Z \exp u(1)$, or 
\item
$N / \Z \exp u(1) \times \Z \exp u(\rho)$ with $\rho \in D_1$. 
\end{enumerate}
For $\rho, \rho' \in D_1$, two Lie groups 
$N / \Z v(1) \times \Z \exp u(\rho)$ and 
$N / \Z v(1) \times \Z \exp u(\rho')$ are isomorphic if and only if 
\begin{itemize}
\item[{\rm (i)}] 
$\{ \rho, \rho' \} = \{ ie^{i \theta}, ie^{-i \theta} \}$ 
for $0 \leq \theta \leq \pi/6$, or 
\item[{\rm (ii)}] 
$\{ \rho, \rho' \} = \{ e^{i \theta}-1, -e^{-i \theta}+1 \}$ 
for $0 < \theta \leq \pi/3$. 
\end{itemize}
\end{Thm}

\begin{proof}
Theorem \ref{maintheorem} states that 
the moduli space of a connected Lie group which is locally isomorphic to $N$ 
is bijective to 
\[ 
\coprod_{r=0}^{2} \ ( \SL_r(\Z) \backslash \SL_r(\R) / \SO_r(\R) ) . 
\] 
Note that 
$\SL_0(\Z) \backslash \SL_0(\R) / \SO_0(\R) = \{ \mathrm{pt} \} = \{ N \}$, 
which corresponds to simply-connected one. 
Furthermore, 
$\SL_1(\Z) \backslash \SL_1(\R) / \SO_1(\R) = \{ \mathrm{pt} \}$, 
which means that $N / \Z u(\rho)$ are isomorphic to each other for all $\rho \in \C$. 
Therefore, in this case, it is isomorphic to $N / \Z u(1)$. 

Now we have to study 
$\SL_2(\Z) \backslash \SL_2(\R) / \SO_2(\R)$. 
It is known that 
$\SL_2(\R) / \SO_2(\R) = \R^{+} \cdot \SL_2(\R) / \CO^+_2(\R)$ 
is bijective to the real hyperbolic plane (or the upper half complex plane) 
\[ 
\R H^2 := \{ z \in \C \ ; \ \im z > 0 \} . 
\] 
Note that the bijective correspondence is induced from 
the following linear fractional transformation: 
\begin{eqnarray*}
\label{ichiji}
p : \R^{+} \cdot \SL_2(\R) \rightarrow \R H^2 : 
\left[ \begin{array}{cc}
a & b \\ c & d 
\end{array} \right] \mapsto (ai+b)/(ci+d) . 
\end{eqnarray*}
Furthermore, note that the action of $\SL_2(\Z)$ on 
$\R^{+} \cdot \SL_2(\R) / \CO^+_2(\R)$ induced from the left action on 
$\R^{+} \cdot \SL_2(\R)$ is equivalent to 
the action of $\SL_2(\Z)$ on $\R H^2$ 
via the linear fractional transformation. 

It is well known (see \cite{Mumford}, \cite{Ume}, for example) that 
representative of the orbit under the action of $\SL_2(\Z)$ on $\R H^2$ 
can be taken as 
\[ 
D_2 := \{ z \in \R H^2 ; \ |z| \geq 1, \ -1/2 \leq \re z \leq 1/2 \} . 
\] 
Note that $z,z' \in D_2$ are in the same $\SL_2(\Z)$-orbit if and only if 
\begin{itemize}
\item[{\rm (i)'}] 
$\{ z, z' \} = \{ ie^{i \theta}, ie^{-i \theta} \}$ 
for $0 \leq \theta \leq \pi/6$, or 
\item[{\rm (ii)'}] 
$\{ z, z' \} = \{ 1/2 + iy, -1/2 + iy \}$ for $y \geq \sqrt{3}/2$. 
\end{itemize}

We translate these results into matrices, 
via the linear fractional transformation $p$. 
Direct calculations show that 
\[
p \left( \left[ 
\begin{array}{cc}
1 & 0 \\ r \cos \theta & r \sin \theta 
\end{array}
\right] \right) = (1/r) e^{i \theta} . 
\]
Therefore the $\SL_2(\Z)$-orbit in $\R^{+} \cdot \SL_2(\R) / \CO^+_2(\R)$ 
have to meet $p(D_3)$, where 
$D_3 \subset \R^{+} \cdot \SL_2(\R)$ is defined by 
\[ 
D_3 := \left\{ \left[ 
\begin{array}{cc}
1 & 0 \\ r \cos \theta & r \sin \theta 
\end{array}
\right] ; \ 
0 < \theta < \pi, \ 0 < r \leq 1, \ 
-1/2 \leq (1/r) \cos \theta \leq 1/2 
\right\} . 
\] 
The above matrix represents the lattice $\Z u(1) + \Z u(re^{i \theta})$ in $\LZ$, 
and hence, represents the discrete subgroup 
$\Z \exp u(1) \times \Z \exp u(re^{i \theta})$ in $Z(N)$. 
The domain $D_1$ can be obtained by the expression of $D_3$, 
and the equivalence conditions (i) and (ii) come from 
the conditions (i)' and (ii)', respectively. 
\end{proof}

\bibliographystyle{amsplain}

\end{document}